\documentclass{article}%
\usepackage{amsmath}
\usepackage{amsfonts}
\usepackage{amssymb}
\usepackage{graphicx}
\usepackage{rotating}
\usepackage[ruled,vlined,linesnumbered,noresetcount]{algorithm2e}
\usepackage{hyperref}
\usepackage{colortbl}
\usepackage[round]{natbib}
\usepackage{mathabx}
\providecommand{\U}[1]{\protect\rule{.1in}{.1in}}
\newtheorem{theorem}{Theorem}

\newtheorem{proposition}[theorem]{Proposition}

\newcommand{\M}{EVWS-TP}

\newenvironment{proof}[1][Proof]{\noindent\textbf{#1.} }{\ \rule{0.5em}{0.5em}}
\begin{document}
\title{Trip Optimization methods for Electric Vehicles Supported by Wireless}
\author{Laura Ant\'on$^{1}$, Marina Leal$^{2}$, Jos\'{e} Luis Sainz-Pardo$^{3}$, \\{\small Centro de Investigaci\'{o}n Operativa, Universidad Miguel
		Hern\'{a}ndez de Elche,}\\{\small Avda. Universidad s/n, 03202, Elche (Alicante), Spain}\\$^{1}${\small l.anton@umh.es}, $^{2}${\small m.leal@umh.es}, $^{3}${\small jose.sainz-pardo@umh.es}\\}
\date{}
\maketitle

\begin{abstract}
	The ecological transition and the reduction of pollution are reasons for the use of green logistic and non-polluting vehicles. In this sense, electric vehicles are an alternative and in recent years many technological advances have been developed by the automotive industry in order to overcome the inherent difficulties in this type of vehicles. In parallel, the classic models and methods for optimizin vehicle routes have been adapted to obtain routings of electric vehicles considering their particular characteristics such as battery autonomy, shortage or recharging stations and elevated time of charge. In this paper, it is studied the optimization of routes taking into account the technological advance due to the possibility of attaching mobile energy diffusers to electric vehicles for recharging them. In this sense, it is developed a model that improves other existing models and, finally, it is also exposed an algorithm to efficiently solve large instances as our computational experience indicates.
\end{abstract}

\bigskip


\bigskip

\textbf{Keywords: }Electric Vehicle Supported by Wireless, trip optimization, Electric Vehicle Routing Problem, constrained shortest path

\section{Introduction}

Nowadays ecological transition promotes the use of non-polluting vehicles as well as green logistics. Government agencies and industries are aware of the need to reduce the environmental impact of vehicles, means of transport and logistic activities. In this sense, several technological advances in electric vehicles (EV) have been developed since EVs are among the cleanest means of transport. These advances offer new possibilities and, at the same time, require new challenges in logistics planning. A new advance, still in development, is the on-the-go charging of vehicles wirelessly. Like the emergence of the electric vehicle and the problem, generally related to its relative autonomy as well as the scarcity of recharging points, it has originated and continues to originate several research works; the possibility of recharging vehicles via wireless brings new challenges.\\

Many classic problems have been adapted to the new problems and challenges posed by electric vehicles. In this way, the Vehicle Routing Problem (VRP) seeks to minimize the total transportation costs of visiting a set of customers for a fleet starting and ending at the depot. Dantzig and Ramsey introduced the problem in 1959 and it has resulted several works such as these of \cite{vrp1}, \cite{vrp2} and \cite{vrp3}. On the other hand, they can be cited variants that take into account usual constraints, then \cite{vrp4} introduced the load capacity of the vehicles, and \cite{vrp5} added customer time windows. Others works can be enumerated as the proposals of \cite{vrp4} about the Vehicle Routing Problem with Backhauls, \cite{vrp8} about the Dynamic Vehicle Routing Problem, etc.\\

Regarding Vehicle Routing Problems, VRP and many of its variants have been transformed into Electric Vehicle Routing Problems (EVRP) taking into account the particularities of these vehicles like autonomy, scarcity of recharging stations, reduction of gas emissions, etc. Just to name a few: \cite{ev1}, \cite{ev2} and \cite{laporte}. Of particular interest is the work of \cite{revisited} that transforms the routing of EV as a Constrained Shortest Path Problem and solve it by Dynamic Programming based on labelling nodes.\\

Finally, our paper is inspired by the work of \cite{original} which studies the incorporation of vehicles that transfer electrical energy to others wirelessly. In this work, the authors propose a model and a heuristic algorithm. Our aim is to improve his research in several aspects.\\




\section{Problem description and formulation}

The models worked in this paper constitute different improvements, even adjustments, over the model introduced in \cite{original}. Therefore and first of all, the cited model is reproduced and analysed. This basically tries to obtain optimal routes for multiple electric vehicles, each of them with an origin and a destination, minimizing the sum of their trip durations. Throughout their journeys, the electric vehicles can be recharged both by mobile energy disseminators (MEDs) and by static charge stations (SCSs). If a MED is attached to an EV this can partially be charged, but they are fully recharged when using a SCS. The problem is defined by a directed graph $G=(N,A)$, not necessarily completed. On the one hand, $N=\{1, ..., n\}$ is the set of vertices that represent transit nodes, SCSs or attachment points of MEDs. On the other hand, $A=\{(i, j) | i, j \in V\}$ defines the set of arcs that represent links between vertices. Other sets are $S \subseteq N$, the set of SCS locations; $M \subseteq N$, the set of MED attachment locations and $K$ the set of EVs. Like in similar routing problems the authors have extended the sets $S$ and $M$ to the sets $S'$ and $M'$ by adding node copies in order to accommodate the possibility of multiple visits to recharging stations and attachment points, later this aspect would be discussed. Regarding the parameters, each arc is doubly weighted by $dt_{ij}$, the driving time between nodes $i$ and $j$, as well as by $c_{ij}$, the energy consumed to traverse the arc $(i,j)$. Some nodes are associated with the parameters $wt_i$, the expected waiting time before using a SCS or MED at location $i$ and $ct_i$, the expected expected charge time to use the SCS $i$. Other parameters are the induced energy factor $\rho$ and $Q^k$, the battery capacity for each EV $k$ and the nodes origin and destination for each vehicle $k$, respectively $s^k$ and $t^k$. Finally, four families of binary decision variables were contemplated: $x_{ij}^k$ which takes a value of $1$ if the vehicle $k$ traverse the arc $(i,j)$, $y_{ij}^k$ which takes a value of $1$ if the vehicle $k$ is receiving energy from a MED when traverse the arc $(i,j)$, $z_i^k$ which takes a value of $1$ if the vehicle $k$ is charged at SCS $i$ and $q_i^k$ which takes a value of $1$ if a MED is attached to the vehicle $k$ at location $i$. Besides, variables $\epsilon_i^k$ save the energy level of each vehicle $k$ when each node $i$ is leaved. These sets, parameters and decision variables are summarized in Table \ref{notacion}. On the other side, routing problems are based on graphs. A graph $G=(V,E)$ consists of a vertex set $V$ and an edge set $E$. Let $N(v) = \{u \in V | (v, u) \in E\}$ and $N[v] = N(v) \cup \{v\}$ be the open and close neighbourhood of $v \in V$, respectively. The degree of a vertex $v \in V$ is $|N(v)|$ and $\delta(G)$ is the minimum degree of a vertex in $G$.\\

Given the above definitions, the mathematical model introduced in \cite{original} was formulated as follows:

\[ \min\ \sum_{k \in K} ( \sum_{(i,j) \in A, i \neq j} (dt_{ij}x_{ij}^k)  + \sum_{i \in S \cup S'} (ct_i + wt_i) z_i^k + \sum_{i \in M \cup M'} (wt_iq_i^k) )  \] \\
\hspace*{1.5cm} s.t.
\begin{eqnarray}
	\sum_{j \in N} x_{ij}^k - \sum_{j \in N}  x_{ji}^k = \left\{
	\begin{array}{rl}
		1 \textrm{, if }i=s^k;\\
		-1\textrm{, if }i=t^k;\\
		0 \textrm{, otherwise,}
	\end{array}
	\right. & & \forall i \in N, \forall k \in K \label{fObj1} \\
	x_{ij} - y_{ij} \geq 0, & & \forall k \in K, \forall j \in N, \forall i \in N, i \neq j \label{flujo1}\\
	\epsilon_j^k \leq \epsilon_i^k - c_{ij} x_{ij}^k + \rho dt_{ij} y_{ij}^k x_{ij}^k + Q^k (1-x_{ij}^k), &  & \forall k \in K, \forall j \in N, \forall i \in N, i \neq j \label{enganche}\\
	\epsilon_i^k \leq Q^k, & & \forall k \in K, \forall i \in N \label{energia1} \\
	\epsilon_i^k = Q^k z_i^k, & & \forall k \in K, \forall i \in S \cup S' \label{energia2} \\
	\epsilon_i^k \geq c_{ij}, & & \forall k \in K, \forall i \in N, \forall j \in S \cup S' \cup M \cup M', i \neq j \label{energia3} \\	
	x_{ij}^k, y_{ij}^k \in \{0,1\} & &  \forall k \in K, \forall i \in N, \forall j \in N,  i \neq j \label{declaracion1} \\	
	z_{ij}^k \in \{0,1\} & &  \forall k \in K, \forall i \in S \cup S' \label{declaracion2} \\
	q_i^k \in \{0,1\} & & \forall k \in K, \forall i \in M \cup M' \label{declaracion3} \\
	\epsilon_i^k \geq 0, & & \forall k \in K, \forall i \in N \label{declaracion4} 
\end{eqnarray}

The objective function \eqref{fObj1} minimizes the sum of total time of the trips. Regarding the constraints, since some of them will be discussed later, they are briefly commented in the following. Constraints \eqref{flujo1} ensure source, destination and conservation flows. Constraints \eqref{enganche} enable the variables $y_{ij}^k$ to take a value of $1$, only if the arc $(i,j)$ is traversed by the vehicle $k$. Constraints \eqref{energia1}-\eqref{energia3} try to set the remaining energy for each vehicle at each node. Finally, \eqref{declaracion1}-\eqref{declaracion4} define the domains of the decision variables. \\

\begin{proposition} \label{separated}
	\cite{original} model can be decomposed into $|K|$ blocks of variables completely separate.
\end{proposition}

\begin{proof}
Trivial. Just note that every decision variable has a superindex $k$ so these conveniently grouped can be structured in $|K|$ different blocks, each of them related to an EV, or what is the same, to a trip.
\end{proof} \\

Given that obtaining the solutions of $|K|$ separate problems is equivalent but more efficient and simple than solving the big problem, then the original model will be transformed into a model related only to one EV and as a consequence, from now the superindices $k$ will be removed. Therefore the solution of the complete problem will be given by the addition of the solutions obtained for the individual trips.

\begin{proposition} \label{Sprima}
The optimal solution of \cite{original} model can be obtained by $S'=\emptyset$.
\end{proposition}

\begin{proof}
Since EV are fully charged by SCSs, a walk returning to the same SCS does not provide any advantage in terms of energy in order to improve the resources to get to the destination, but it does imply a greater duration of the trip thus worsening the objective value. 
\end{proof}\\

\begin{figure}[!h]
	\begin{center}
		\includegraphics[width=0.9\textwidth]{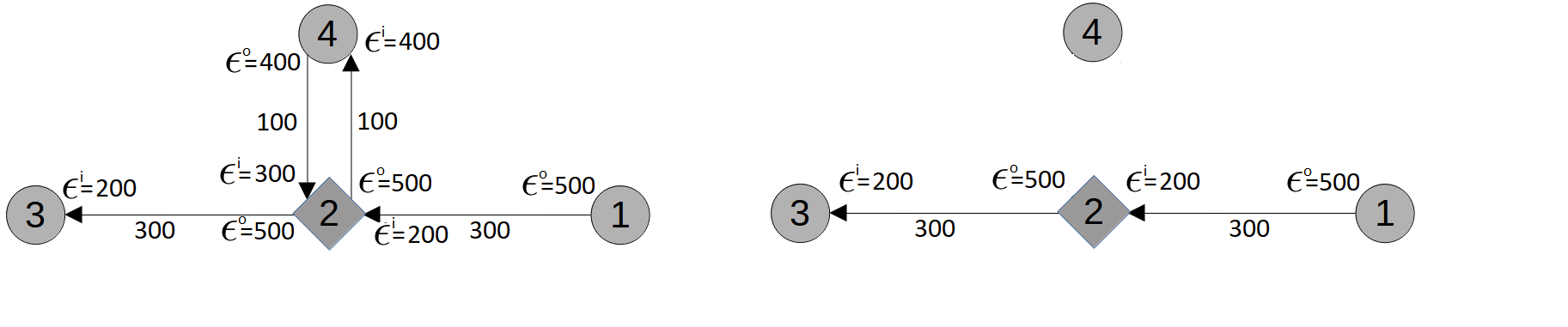}
		\caption{Example of inefficiency of returning to the same SCS \label{SPrima}}
	\end{center}
\end{figure}

Figure \ref{SPrima} illustrates this proposition. The origin node of the trip is $1$ and the destination $3$ while $4$ is a transit node and $2$ a SCS. The values on the arcs represent the cost of traversing them in terms of energy and the values $\epsilon^0$ and $\epsilon^i$ respectively denote the energy of departure and arrival at each node. Suppose the energy capacity of the vehicle is $Q=500$, initially the EV is fully charged and departs with a remaining energy of $\epsilon^o = 500$. It can be seen that the graph on the left represents the non-optimal path $(1,2) - (2,4) - (4,2) - (2,3)$. Since the EV is fully charged at node $2$ then all routes that involve returning to SCS $2$ do not provide more power to the vehicle and lead to a longer delay, and even higher energy consumption, to reach the destination. Therefore, it is more efficient to continue to the destination without going through the SCS $2$ twice, as it is shown in the graph to the right. In this case, given that $Q=500$ and $c_{23}=300$, the destination can be reached directly, but in the case where the EV did not have enough energy for it, it would be necessary to continue the trip trying to reach other SCS or MED before to get to node $3$.

\begin{proposition} \label{separated}
	If the subset $N \setminus (M \cup S)$ were not extended with copies then the optimal solution of \cite{original} model could not be obtained.
\end{proposition}

\begin{figure}[!h]
	\begin{center}
		\includegraphics[width=0.4\textwidth]{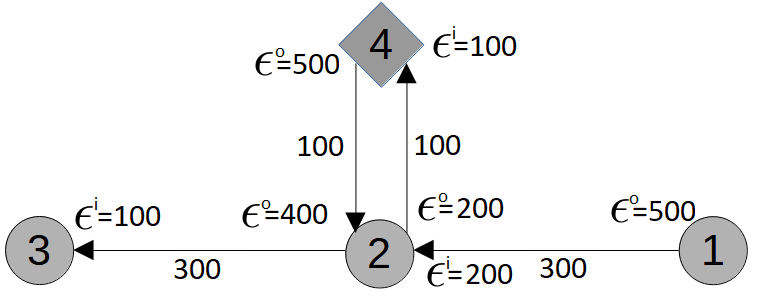}
		\caption{Counter-example of solution not found in case of not extending the set $N$.\label{NPrima}}
	\end{center}
\end{figure}

\begin{proof}
The subset $N \setminus (M \cup S)$ just refers to transit nodes, that is, nodes that do not have SCS and are not location for MED connections. Figure \ref{NPrima} illustrates a case in which it is more advantageous to traverse the transit node $2$ twice. When the EV reaches $2$ from $1$ it has not enough power to get to node $3$, since $200<300$. But it can get to SCS $4$, fully recharge its battery, and then go back to node $3$ with an energy level of $400$, which is enough power to get to node $3$. Since it may be more efficient to go through transit nodes several times, then it is also necessary to extend the subset $N \setminus (M \cup S)$.
\end{proof}

At this point, it is introduced the Electric Vehicle Wireless Supported Trip Problem (\M) which is based on the \cite{original} model. In addition to improving and adjusting some of its sets and constraints, other aspects have been taken into account. Firstly, a MED cannot provide unlimited energy to an EV. Secondly, a reduction in circulation speed  while a MED keeps on attaching to an EV has also been contemplated which seems more realistic. So, the Electric Vehicle Wireless Supported Trip Problem is formulated as:\\

\[ \mbox{(\M)} \min\ \sum_{i \in N'A} \sum_{j \in N'} d_{ij} x_{ij}  + \sum_{i \in N'A} \sum_{j \in N'} \sum_{m \in M'}\frac{\delta d_{ij}}{1-\delta}  y_{ij}^m +  \sum_{i \in S}  \frac{t_i}{Q} \beta_i + \sum_{i \in S} w_i z_i + \sum_{i \in M'} w_i \hspace{0.1cm} q_i \] \label{fObj2}\\
\hspace*{1.5cm} s.t.
\begin{eqnarray}
	\sum_{j \in N'} x_{ij} - \sum_{j \in N'}  x_{ji} = \left\{	
	\begin{array}{rl}
		1 \textrm{, if }i=s;\\
		-1\textrm{, if }i=t;\\
		0 \textrm{, otherwise,}
	\end{array}
	\right. & & \forall i \in N' \label{flujo2} \\
	\sum_{j \in N'} y_{ij}^m - \sum_{j \in N'}  y_{ji}^m \left\{
	\begin{array}{rl}
		\leq q_i \textrm{, if }i=m;\\
		= -a_i^m\textrm{, otherwise}
	\end{array}
	\right. & & \forall i \in N', \forall m \in M'  \label{MEDflujo} \\
	\sum_{i \in N': i \neq m} t_i^m=q_m & &\forall m \in M' \label{MEDflujo2} \\
	\gamma_{ij}^m \leq \frac{\rho}{1-\delta} d_{ij} y_{ij}^m & & \forall i \in N', \forall j \in N', \forall m \in M' \label{MEDflujo3}\\
	\sum_{i\in N'}	\sum_{j\in N'} \gamma_{ij}^m \leq W& & \forall m \in M' \label{MEDflujo4}\\
	\sum_{i \in S} q_i \leq P & & \label{Prestriccion}\\	
	\sum_{j\in N'} \sum_{m \in M'} y_{ij}^m \leq Q - \epsilon_i & & \forall i \in N' \label{MEDflujo5}\\
	z_i \leq \beta_i \leq Q z_i & & \forall i \in S \label{Ecargada} \\	
	\beta_i = 0 & & \forall i \in N' \setminus S \label{Beta0}\\
	\epsilon_j \leq \epsilon_i - c_{ij}  x_{ij} + \sum_{m \in M'} \gamma_{ij}^m + \beta_j + (1-x_{ij}) Q & & \forall i \in N', \forall j \in N' \label{energia1} \\
	\epsilon_i - c_{ij} x_{ij}  + \sum_{m \in M'} \gamma_{ij}^m + (1-x_{ij}) Q \geq 0\label{energia2}  & & \forall i \in N', \forall j \in N'\\
	\epsilon_i \leq Q \label{energia3} & & \forall i \in N'\\
	\epsilon_i \geq Q z_i \label{energia4} & & \forall i \in S'\\
	\epsilon_0 = Q \label{energia5} & & \\
	z_i \leq \sum_{j \in N'} x_{ji}\label{opcional1}  & & \forall i \in S\\
	z_i \leq \sum_{j \in N'} x_{ij} \label{opcional2}  & & \forall i \in S\\	
	x_{ij}^m, y_{ij}^m \in \{0,1\} & & \forall i \in N, \forall j \in N,  i \neq j, \forall m \in M' \label{declaracion1} \\	
	z_{i} \in \{0,1\} & & \forall i \in S' \label{declaracion2} \\
	q_i \in \{0,1\} & & \forall i \in  M' \label{declaracion3} \\
	\epsilon_i, \gamma_i \geq 0, & & \forall i \in N' \label{declaracion4} \\
	a_i^m \in \{0,1\} & & \forall m \in M', \forall i \in N' 
\end{eqnarray}

\begin{table}[h]
	\centering
	\caption{Mathematical notation of \M}
	\label{notacion}
	\begin{tabular}
		[c]{ll}
		\hline
		\textbf{Sets} & \\
		$N$ & Set of nodes \\
		$S \subseteq N$ & Set of SCSs \\
		$M \subseteq N$ & Set of locations where a MED can be connected\\	
		$N'$ & Set of nodes with their copies \\
		$S' \subseteq N'$ & Set of SCSs with their copies \\		
		$M' \subseteq N'$ & Set of locations where a MED can be connected with their copies\\
		\textbf{Parameters} & \\
		$d_{ij}$ & Driving time to traverse the arc $(i,j)$ \\
		$c_{ij}$ & Energy to traverse the arc $(i,j)$\\
		$w_i$ & Expected waiting time before using a SCS or MED at location $i$\\
		$t_i$ & Full charge time at SCS $i$\\
		$b_i$ & Battery full charge time at MED $i$\\		
		$\rho$ & Factor of induced energy charged from a MED (kWh per time unit)\\
		$\delta$ & Factor of speed reduction when a MED is connected\\
		$Q$ & Battery capacity of EV\\
		$s$ & source node for the $k$ vehicle trip\\ 
		$t$ & destination node for the $k$ vehicle trip\\
		$W$ & Limit of induced energy from a MED\\
		\textbf{Decision variables} & \\
		$x_{ij}$ & $1$ if EV departs from node $i$ and arrives at node $j$, $0$ otherwise\\
		$y_{ij}$ & $1$ if EV traverses arc $(i,j)$ attached to a MED, $0$ otherwise\\
		$z_i$ & $1$ if EV is charged at SCS $i \in S'$, $0$ otherwise\\
		$q_i$ & $1$ if a MED is attached to the EV at location $i \in M'$, $0$ otherwise\\
		$a^m_i$ & $1$ if MED $m$ is disconnected at node $i$, $0$ otherwise\\				
		$\epsilon_i$ & battery SoC $k$ when node $i$ is leaved\\
		$\beta_i$ & battery charged at node $i$ (if node $i$ is not a SCS then $\beta_i=0$)\\	
		$\gamma_{ij}^m$ & energy induced by MED $m$ traversing $(i,j)$\\
		\hline
	\end{tabular}
\end{table}

The sets, parameters and decision variables of the problem are presented in Table \ref{notacion}. Note the model presents a slight overuse of notation, especially if the graph is not completed but also contemplating arcs with the same origin and destination. As it can be seen sets $M$ and $N \setminus M$ have been extended with its corresponding copies to contemplate multiple visits. Given that an efficient policy is not to be connected to a MED when the battery is fully charged then $S \cup M = \emptyset$ and $S \cup M' = \emptyset$. The objective function continues minimizing the total duration of the trip. The summation $\sum_{(i,j) \in A} dt_{ij} \left(x_{ij}  + \sum_{m \in M'}\frac{\delta }{1-\delta}  y_{ij}^m \right)$ represent the total time of traversing the arcs of the trip. On the one hand, $\sum_{(i,j) \in A} d_{ij} \left( x_{ij} - \sum_{m \in M'} y_{ij}^m \right)$ is the time of traversing the arcs without being connected to a MED. On the other hand, $\sum_{m \in M'}\frac{d_{ij}}{1-\delta}  y_{ij}^m$ is the time of traversing the arcs being connected to a MED considering that the EV speed have been reduced a $100 \hspace{0.05cm} \delta \hspace{0.05cm} \%$ . Since, $\sum_{(i,j) \in A} d_{ij} \left( x_{ij} - \sum_{m \in M'} y_{ij} \right) + \sum_{(i,j) \in A}  \sum_{m \in M'}\frac{d_{ij}}{1-\delta}  y_{ij}^m  = \sum_{(i,j) \in A} d_{ij} x_{ij}  + \sum_{(i,j) \in A} \sum_{m \in M'}\frac{\delta d_{ij}}{1-\delta} y_{ij}^m$ this represent the total time of traversing the arcs. The rest of terms of the objective function refer to the total charging time plus the total waiting time.\\

Constraints \eqref{flujo2} ensure the flow between origin $s$ and destination $t$. Another group of flow constraints \eqref{MEDflujo}-\eqref{MEDflujo5} have been introduced in order to limit the energy induced by each MED. Besides, if MED $m$ is used, that is $q_m=1$, then it is also obtained its disconnection node $i$ such that $t_i^m=1$. If EV traverses $(i,j)$ attached to a MED then the induced energy $\gamma_{ij}^m$ will be at most the time employed to traversing this arc, $\frac{d_{ij}}{1-\delta}$ multiplied by the factor of induced energy $\rho$ as it is expressed by $\eqref{MEDflujo3}$. Given that a MED cannot provide unlimited energy, the total induced energy of each MED $m$ is constrained by $W$ at \eqref{MEDflujo4}. Besides, \eqref{Prestriccion} prevents to overcome $P$ connections to different MEDs. \eqref{MEDflujo5} imposes the policy of not being connected to a MED if the battery is full, that is, $\epsilon_i=Q$. Finally, the family of constraints \eqref{Ecargada} indicates that EV will be recharged by a SCS only if this provide it some energy but without overcoming its battery capacity. Given that it is necessary to define $\beta_i$ variables for all $i \in N'$ they are set to $0$ in \eqref{Beta0} for the nodes that are not SCS.\\

Regarding the energy flow this is expressed in terms of energy level when the EV leaves a node. If $\epsilon_i$ is the energy level when EV departs from node $i$ and traverses arc $(i,j)$, then the energy level $\epsilon_j$ at the moment EV departs from node $j$ will be given by expression \eqref{Ecenergia}.

\begin{equation} \epsilon_j = \epsilon_i - c_{ij}  + \sum_{m \in M'} \gamma_{ij}^m + \beta_j  \label{Ecenergia} \end{equation}

Being $c_{ij}$ the consumed energy, $\sum_{m \in M'} \gamma_{ij}^m$ the induced energy and $\beta_j$ the charged energy at node $j$. This relation is expressed by $\eqref{energia1}$, but the constraint is disabled in case of not traversing arc $(i,j)$ adding the term $(1-x_{ij}) Q$. \eqref{energia2} imposes the EV has to have enough energy before traversing the arc $(i,j)$. Finally, the energy level at each node have to be under the battery capacity \eqref{energia3}, but in case of charging the EV have to fully charged \eqref{energia4} and the EV departs from the origin node fully charged \eqref{energia5}.\\

Constraints \eqref{opcional1} and \eqref{opcional2} are not compulsories but they help to find the optimal solution. On the one hand, a SCS $i$ is not used if the EV does not arrive to it. On the other hand, a SCS $i$ is not used if the EV does not depart from it. \\


\begin{proposition} 
	\M\ is NP-hard.	
\end{proposition}

\begin{proof}
	Given that instances where $S=\emptyset$ and $M=\emptyset$ are equivalent to constrained shortest path problems, \M\ generalizes the constrained shortest path problem which is NP-hard.
\end{proof}

\section{Optimal strategies and dominated solutions}

As the solution method to the \M\, an exact algorithm based on Dynamic Programming (DP) is proposed. This algorithm like most of the DP methods used for routing are based on node labelling. But like branch and cut methods, it also includes pruning the pool of solutions to be explored according to dominance criteria and bounds obtained from the solutions already found; so, it is possible to obtain the optimal solution before labelling all the nodes. Pruning by dominance relationships have been also developed in works like in \cite{constrained} for efficiently solving the constrained shortest path but, although there are works like in \cite{revisited} that transforms energy-optimal routing without recharging problems as special cases of the constrained shortest path introducing even DP algorithms to solve it, the /M/ cannot be transformed in a constrained shortest path due to the possibility of effectuating recharges, in any case, it may be transformed in a constrained shortest path with dynamic and negative costs. There are also works, employ dominance criteria to eliminate stations in a preprocessing part like \cite{laporte}. But, as novelty, pruning dominance relationships for electric vehicle routes is applied in this work.\\

The dominance here exposed is based on walks. Suppose two solutions both traversing a node $w$. In the first solution $w$ is reached by a faster walk and besides the EV also has more available resources: energy and disposable MEDs, obviously the second solution can be improved going from origin to $w$ by the walk of the first one. Note that has not happen just for the fact of arriving faster to $w$, or even for arriving faster and with more energy, it is also necessary to arrive with less MEDs used or in case of arriving with the same number of MEDs but in case of are using one in both solutions when the EV arrives to $w$, then the MED of the first solution has to have less transferred energy.\\ 

Formally, the relationship of dominance can be exposed as it follows. For a solution $\mathcal{S}$ that traverses node $w$, let $\mathcal{C}^\mathcal{S}_w=\{(s, v_1), ..., (v_{p}, w)\}$ be the partial walk from origin $s$ to node $w$, $\nu^\mathcal{S}_w=\{s, v_1, ..., v_{p}, w\}$ the set of nodes traversed by this partial walk and $\epsilon^\mathcal{S}_i$, $q^\mathcal{S}_i$ the variables about battery levels and MED attachments. Be $t^\mathcal{S}_w$ the time of walking $\mathcal{C}^\mathcal{S}_w$ including waiting and charging times and $a_{w}^\mathcal{S}=m \in M': (\exists j\in N': y^m_{wj}=1)$, that is the MED connected to the EV when it departs from $w$. Note the total induced energy by a MED $m \in M'$ can be computed by $g^{\mathcal{C}^\mathcal{S}_w}_m=\sum_{(i,j) \in \mathcal{C}^{\mathcal{S}}_w} \gamma_{ij}^m$ if $m$ traverses $w$; $g^{\mathcal{C}^\mathcal{S}_w}_m=0$ otherwise. Then, given two solutions $\mathcal{S}$ and $\mathcal{S'}$  both traversing $w$, it is said that walk $\mathcal{C}^\mathcal{S}_w$ dominates walk  $\mathcal{C}^\mathcal{S'}_w$, or $\mathcal{C}^\mathcal{S}_w \succeq \mathcal{C}^\mathcal{S'}_w$ ,if it accomplish $t^{\mathcal{S}_w}_w \leq t^{\mathcal{S'}_w}_w$, $\epsilon^{\mathcal{S}}_w - \beta^{\mathcal{S}}_w \geq \epsilon^{\mathcal{S'}}_w - \beta^{\mathcal{S'}}_w$ and one of the following conditions:\\

\begin{itemize}
	\item $\displaystyle \sum_{i \in  \nu^\mathcal{S}_w  \setminus w} q^\mathcal{S}_i \leq \sum_{i \in  \nu^\mathcal{S'}_w \setminus w} q^\mathcal{S'}_i$
	\item $\displaystyle \sum_{i \in \nu^\mathcal{S}_w  \setminus w} q^\mathcal{S}_i = \sum_{i \in \nu^\mathcal{S'}_w  \setminus w} q^\mathcal{S'}_i$, and $g^{\mathcal{C}^\mathcal{S}_w}_{a_{w}^\mathcal{S}} \leq g^{\mathcal{C}^\mathcal{S'}_w}_{a_{w}^\mathcal{S'}}$
\end{itemize}

Table \ref{sVariables} summarizes the sets and state variables used in the developed algorithm. This DP build two sets of solutions. On the one hand, the set of partial solutions, each one of them not completely explored. On the other hand, the set of dominant solutions which purpose is to cast partial solutions away without ending its exploration given that are dominated. Basically it is an algorithm for traversing the graph, like Depth First Search (DFS) or Breadth First Search (BFS), although the selection of the node to be explored is different from both. Through which they are compiled nodes to be explored and dominant solutions.\\

\begin{table}[h]
	\centering
	\caption{Sets and state variables}
	\label{sVariables}
	\begin{tabular}
		[c]{ll}
		\hline
		\textbf{Sets} & \\
		$\mathcal{A}$ & Set of partial solutions pending to be explored\\
		$\mathcal{B}$ & Set of dominating partial solutions\\
		\textbf{State variables} & \\
		$\nu$ node to \\
		\textbf{Associated state variables} & \\
		$\tau$ & time spent to arrive to the current node\\		
		$\mu$ & Number of used MEDs \\
		$\epsilon$ & Level of energy at the moment of arriving to the node\\				
		$\omega$ & 1, if a MED is attached to the EV when arrives to the current node; 0, otherwise\\		
		$\eta$ & Transferred energy by the current attached MED \\		
		
		\hline
	\end{tabular}
\end{table}

The point of the algorithm is to simulate several feasible walks in order to arrive from $s$ to $t$. Given that at each step one partial solution or walk belonging to $\mathcal{A}$ is explored by incorporating to $\mathcal{A}$ all the decisions that can be made from its current situation. So, each item belonging to $\mathcal{A}$ is made up by the last node reached $\nu$, the spent time to reach $\nu$, the energy level of the EV $\epsilon$, $\omega$ or if a MED is attached to the EV when it arrives to $\nu$  and also $\eta$, the energy transferred to the EV by its current attached MED. Then each item $a \in \mathcal{A}$ is defined by $(\nu, \tau, \mu, \epsilon, \omega, \eta)$. \\

Dynamic programming (DP) is a methodology commonly applied to problems that can be divisible into stages in each of which decisions are made, with the aim of minimizing a cost function. A deterministic state usually is expressed as $x_{k+1}=f(x_k, u_k)$ where $u_0, ..., u_{N-1}$ are the control inputs. Each control input is composed of state variables whose values, often discrete and finite,  reflect the decisions made. They can also introduced associated state variables which are not referred to decisions but values that take other variables in function of the decision made. Regarding the cost function, also named policy in the DP jargon, this commonly is of the form: $J=g_N(x_N) + \sum_{k=0}^{n-1} g_k(x_k, u_k)$. DP is based on the Bellman Optimity Principle. Accordingly, the optimal policy can be obtained by a backward approach that solves the last stage first, after the subproblem involving the last two stages, and so on until the whole problem is finally solved.

In our case, one or more of the following decisions can be made at each node

\begin{itemize}
	\item Get to a neighbour.
	\item Attach a MED.
	\item Fully recharge the battery.
\end{itemize}

The feasibility of each decision depends on the type of the node and the estate variables of the EV. The algorithm proposed also precise two type of sets. The set $\mathcal{U}$ formed by a partial solution pending of to be explored and de set $\mathcal{D}$ formed by arrays of partial solutions that are dominants in the sense of...

The main loop consists in selecting one partial walk belonging to $\mathcal{A}$ and to obtain all the decisions that can be make from the situation of the EV in the last node $\nu$ reached. The feasible decisions obtained will be incorporated to $\mathcal{A}$ as new items to be explored. So, all the neighbours of $\nu$ where it is possible to arrive with the current levels of energy are added to $\mathcal{A}$ by updating the state variable values. In case of arriving to $\nu$ with a MED attached, then it will be added all the neighbours reachable without and with the MED conveniently updated the state variables. If there is a SCS at $\nu$ to charge the EV should also be considered. Regarding the criteria to select the item to be explored these are neither LIFO or FIFO criteria...
 
Initially \eqref{inicio}, it is set $\mathcal{A}=\{a_0=(s, 0, 0,Q, 0, 0) \}$, that is, the EV is at location $s$ full of energy and without a MED attached. Then, the main loop \eqref{cicloIni}-\eqref{cicloFin} consists on exploring all the walks or partial solutions belonging to $\mathcal{A}$. So, at the beginning $a_0$ is explored, incorporating into $\mathcal{A}$ all the feasible walks obtained by contemplating all the possible decisions. All the neighbours of $s$ that can be reached with the current level of energy are incorporated to $\mathcal{A}$ with their state variable values of time, energy, etc. updated according to this decision. It will continue selecting another partial walk $a_i$ to be explored. And all the combinations of possible decisions that can be made from the situation of the last node that was reached in this walk will be annotated in order to be explores, for example, all the neighbours that can be reached charging and without energy if this node were a SCS, and so on. Decir at each node it can be made up to three types of decisions... Once a partial walk is explored advancing it then it is erased from $\mathcal{A}$.\\

\begin{algorithm}[h]
	\setcounter{AlgoLine}{0}
	\label{initial}
	\caption{ECVW exploration}
	\ForEach{$v \in V$}
	{
		$d(v) = (+\infty, +\infty, +\infty)$
	}
	$\mathcal{P}_0 = \{O\}$\\
	$\mathcal{A}=\{a_0=(\nu=s, \tau=0, \mu=0, \epsilon=Q, \omega=0, \eta=0, \mathcal{P}_0\}$ \label{inicio}\\
	\While{$\mathcal{A} \neq \emptyset$  \label{cicloIni}}
	{
		$u' = u \in \mathcal{A} : dijkstra(u', D) = min(dijkstra(u, D), \forall u \in U)$\\
		\ForEach{$v \in N(u')$}
		{
			$\mathcal{P}=\mathcal{P}_{u'} \cup v$\\
			$\tau = \tau + d_{\nu,v}$\\
			$\nu = v$\\
			$\epsilon = \epsilon $\\			
			\If{$v \in S$}
			{
			}
			\If{$v \in M'$}
			{
			}
			$d=d_{u'} + d_{u'v}$\\
			$\epsilon = \epsilon_{u'} - c^1_{u'v}$\\
			$\gamma = \gamma_{u'} - c^2_{u'v}$\\
			$u_0=(v, d, \epsilon, \gamma, 0, \mathcal{P})$\\
			\If{$u_0$ is feasible}
			{
				\uIf{$v=D$}
				{
					\If{$d<B$}
					{
						$B=d$\\
						$\mathcal{P}*=\mathcal{P}$\\
					}							
				}
				\uElseIf{$d < B$ \textbf{and} $\mathcal{U} \succeq u_0$ \textbf{and}  $\mathcal{E} \succeq u_0$}
				{
					$\mathcal{A}=\mathcal{A} \setminus \{u \in \mathcal{A}: u_0 \succ u \}$\\								
					$\mathcal{A}=\mathcal{A} \cup u_0$\\								
				}
				$u_1=(v, d, \epsilon, \gamma, 0, \mathcal{P})$\\	
				\If{$\mathcal{U} \succeq u_1$ \textbf{and}  $\mathcal{E} \succeq u_1$}
				{
					$\mathcal{A}= \mathcal{A} \cup u_1$\\								
				}
			}
		}\label{cicloFin}
		$\mathcal{B}=\mathcal{B} \cup u'$\\
		$\mathcal{B}=\mathcal{B} \setminus \{u \in \mathcal{B}: u' \succ u \}$\\			
		$\mathcal{A}=\mathcal{A} \setminus u'$\\		
	}
\end{algorithm}

\section{Computational experience}


We tested the proposed algorithm on 4 instances of 40000 nodes and 100000, 200000, 400000, 600000 and 800000 edges randomly generated. In $10 \%$ of the nodes was randomly setted a recharging station and in $90 \%$ of the nodes was located a MED access point. Finally, a limit of $3$ MEDs was allowed to be used.\\

Table \ref{computacional} shows the results of this computational experience. Its columns represent:

\begin{itemize}
	\item \textit{\#Edges} indicates the number of edges,
	\item  \textit{t} the computational time in seconds employed for solving each instance,
	\item \textit{Energy} the recharged energy,
	\item  \textit{W. Energy} the energy obtained from MEDS,
	\item \textit{\# Recharges} the number of recharges carried out,
	\item \textit{\# MED} the number of MED's employed.
\end{itemize}

\begin{table}[h]
	\centering
	\begin{tabular}{|r|r|r|r|r|r|r|r|r|r|}
		\hline
		\textbf{\#Edges} &  \textbf{t} & \textbf{Kms} & \textbf{Energy} & \textbf{W. Energy} & \textbf{\# Recharges} & \textbf{\# MED} \\ \hline
		 100000      & 15         & 1345         & 0                 & 0                      & 0                 & 0                \\ \hline
		 200000       & 28         & 1257         & 767               & 282                    & 1                 & 3                \\ \hline
		400000        & 56         & 644          & 182               & 182                    & 0                 & 2                \\ \hline
		 600000        & 80         & 337          & 0                 & 0                      & 0                 & 0                \\ \hline
		800000        & 111        & 266          & 0                 & 0                      & 0                 & 0                \\ \hline
	\end{tabular}
	\caption{Computational results} \label{computacional}
\end{table}

\begin{figure}[h]
	\begin{center}
		\includegraphics[width=14cm]{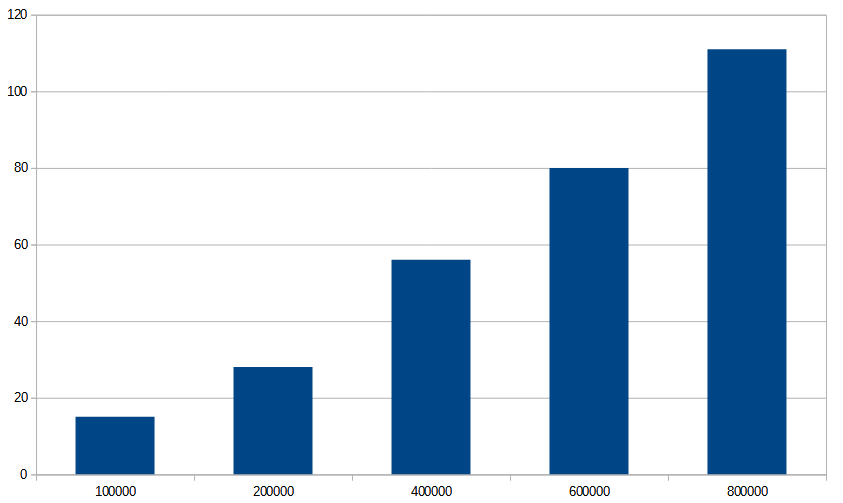}
	\end{center}
	\caption{Time for instance} \label{tiempos}
\end{figure}

Figure \ref{tiempos} represent the time used for solving each instance.  As it can be noted, the times are relatively reduced so the proposed algorithm can be used for solving large instances in an exact way. No instances were solved from the model in less than 3600 s. 

\section{Conclusions}

In this work a model for planning trips of an electric vehicle supported by wireless recharge has been developed. An exact algorithm based on dynamic programming has also been developed. This algorithm takes into account the properties and specific situations about the problem, that have been analysed too. The times for solving the instances generated by the exact algorithm proposed are relatively reduced whilst no instances were solved from the model in less than 3600 s.

\bigskip

\textbf{Acknowledgments}

This work has been partially supported by the Spanish Ministry of Science, Innovation and Universities, project PGC2018-099428-B-I00

\bibliographystyle{apalike}
\bibliography{ev}

\begin{thebibliography}{}

\bibitem[Artmeier et~al., 2010]{revisited}
Artmeier, A., Haselmayr, J., Leucker, M., and Sachenbacher, M. (2010).
\newblock The shortest path problem revisited: Optimal routing for electric
  vehicles.
\newblock In Dillmann, R., Beyerer, J., Hanebeck, U.~D., and Schultz, T.,
  editors, {\em KI 2010: Advances in Artificial Intelligence}, pages 309--316,
  Berlin, Heidelberg. Springer Berlin Heidelberg.

\bibitem[Christofides, 1976]{vrp1}
Christofides, N. (1976).
\newblock The vehicle routing problem.
\newblock {\em Revue française d’automatique, d’informatique et de
  recherche opérationnelle}, 10(v1):55--70.

\bibitem[D.~Pecin and Uchoa, 2017]{vrp4}
D.~Pecin, A.~Pessoa, M.~P. and Uchoa, E. (2017).
\newblock Improved branchcut-and-price for capacitated vehicle routing.
\newblock {\em Mathematical Programming Computation}, 9(1):61--100.

\bibitem[Erdogan and Miller-Hooks, 2017]{ev2}
Erdogan, S. and Miller-Hooks, E. (2017).
\newblock The electric vehicle-routing problem with time windows and recharging
  stations.
\newblock {\em Transportation Research E Logist Transp.}, 48 (1):100--114.

\bibitem[Keskin et~al., 2019]{laporte}
Keskin, M., Laporte, G., and Çatay, B. (2019).
\newblock Electric vehicle routing problem with time-dependent waiting times at
  recharging stations.
\newblock {\em Computers $\&$ Operations Research}, 107:77--94.

\bibitem[Kosmanos et~al., 2018]{original}
Kosmanos, D., Maglaras, L.~A., Mavrovouniotis, M., Moschoyiannis, S., Argyriou,
  A., Maglaras, A., and Janicke, H. (2018).
\newblock Route optimization of electric vehicles based on dynamic wireless
  charging.
\newblock {\em IEEE Access}, 6:42551--42565.

\bibitem[Lozano and Medaglia, 2013]{constrained}
Lozano, L. and Medaglia, A.~L. (2013).
\newblock On an exact method for the constrained shortest path problem.
\newblock {\em Computers $\&$ Operations Research}, 40(1):378--384.

\bibitem[Pillac et~al., 2013]{vrp8}
Pillac, V., Gendreau, M., Guéret, C., and Medaglia, A.~L. (2013).
\newblock A review of dynamic vehicle routing problems.
\newblock {\em European Journal of Operational Research}, 225(1):1--11.

\bibitem[R.~Baldacci and Roberti, 2012]{vrp5}
R.~Baldacci, A.~M. and Roberti, R. (2012).
\newblock Recent exact algorithms for solving the vehicle routing problem under
  capacity and time window constraints.
\newblock {\em European Journal of Operational Research}, 218(1):1--6.

\bibitem[Schneider et~al., 2014]{ev1}
Schneider, M., Stenger, A., and Goeke, D. (2014).
\newblock The electric vehicle-routing problem with time windows and recharging
  stations.
\newblock {\em Transportation Science}, 48 (4):500--520.

\bibitem[Toth and Vigo, 2002]{vrp2}
Toth, P. and Vigo, D. (2002).
\newblock {\em The Vehicle Routing Problem}, volume~9.

\bibitem[Toth and Vigo, 2014]{vrp3}
Toth, P. and Vigo, D. (2014).
\newblock {\em Vehicle Routing: problems, methods and applications}.

\end{thebibliography}

\end{document}